\documentclass [a4paper]{article}
\usepackage{amsmath}
\usepackage{mathrsfs}
\usepackage{amsfonts}
\usepackage {graphicx}
\usepackage{color,xcolor}
\usepackage{amsmath,amssymb}
\input{epsf.sty}
\usepackage{epsfig}

\usepackage{amsthm}
\newtheorem{rem}{Remark}
\newtheorem{theorem}{Theorem}

\def\be{\begin{equation}}
\def\ee{\end{equation}}
\def\br{\begin{eqnarray*}}
\def\er{\end{eqnarray*}}

\begin{document}
\begin{center}
\Large\bf{A theoretical contribution to the fast implementation of null linear discriminant analysis method using random matrix multiplication with scatter matrices
}\\
%
\quad\\

\normalsize~Ting-ting Feng\footnote[1]{School of Mathematics and statistics, Jiangsu Normal University, Xuzhou, 221116, Jiangsu, P.R. China.
Email: {\tt tofengtingting@163.com}. This author is supported by the Postgraduate Innovation Project of Jiangsu Province under grant CXLX13\_968.},
~
Gang Wu\footnote[2]{Corresponding author (G. Wu). Department of Mathematics,
China University of Mining and Technology \& School of Mathematics and statistics, Jiangsu Normal University, Xuzhou, 221116, Jiangsu, P.R. China.
Email: {\tt gangwu76@126.com} and {\tt wugangzy@gmail.com}. This author is
supported by the National Science Foundation of China under grant 11371176, the Natural Science Foundation of Jiangsu Province under grant BK20131126, the 333 Project of Jiangsu Province, and the Talent Introduction Program of China University of Mining and Technology.}

\end{center}

\begin{abstract}
The null linear discriminant analysis method is a competitive approach for dimensionality reduction. The implementation of this method, however, is computationally expensive. Recently, a fast implementation of null linear discriminant analysis method using random matrix multiplication with scatter matrices was proposed. However, if the random matrix is chosen arbitrarily, the orientation matrix may be rank deficient, and some useful discriminant information will be lost.
In this paper, we investigate how to choose the random matrix properly, such that the two criteria of the null LDA method are satisfied theoretically.
We give a necessary and sufficient condition to guarantee full column rank of the orientation matrix. Moreover, the
geometric characterization of the condition is also described.
\\
\mbox{\bf Keywords:} Dimensionality reduction, Linear discriminant analysis (LDA), Null linear discriminant analysis (Null LDA), Small sample size problem.


\end{abstract}

\bigskip

\section{Introduction}
\label{sec1}

\setcounter{equation}{0}

Dimensionality reduction has become an ubiquitous preprocessing step in many applications. In general, its objectives are to remove irrelevant and redundant
data to reduce the computational cost and to improve the quality of data for efficient data-intensive processing tasks such as face recognition and data mining.
Linear discriminant analysis (LDA) is one of the most popular and powerful
dimensionality reduction techniques for classification (Fukunaga, 1990). However, a main disadvantage of LDA is that
the so-called total scatter matrix must be nonsingular. Indeed, in many applications, the scatter matrices
can be singular since the data points are from a very high-dimensional space, and thus usually
the number of the data samples is much smaller than the data dimension. This is the well-known small sample size (SSS) problem or the undersampled
problem (Fukunaga, 1990).

Let $X=[{\bf x}_{1},{\bf x}_{2},\ldots,{\bf x}_{n}]$ be a set of training samples in a $d$-dimensional feature space, and $\Omega=\{\omega_{j}:j=1,2,\ldots,c\}$ be the class labels, with $\omega_{j}$ being the $j$-th class. We denote by $n_{j}$ the number of samples in the $j$-th class, which satisfies $\sum_{j=1}^{c}n_{j}=n$. Let ${\bf \mu}_{j}$ be the centroid of the $j$-th class, and ${\bf \mu}$ be the global centroid of the training data set.
Then we define the within-class scatter matrix
$$
S_{W}=\sum_{j=1}^{c}\sum_{{\bf x}_{i}\in \omega_{j}}({\bf x}_{i}-{\bf \mu}_{j})({\bf x}_{i}-{\bf \mu}_{j})^{T},
$$
and the between-class scatter matrix
$$
S_{B}=\sum_{j=1}^{c}n_{j}({\bf \mu}_{j}-{\bf \mu})({\bf \mu}_{j}-{\bf \mu})^{T}\equiv BB^{T},
$$
where
$B=[\sqrt{n_{1}}({\bf \mu}_{1}-{\bf \mu}),\sqrt{n_{2}}({\bf \mu}_{2}-{\bf \mu}),\ldots,\sqrt{n_{c}}({\bf \mu}_{c}-{\bf \mu})]\in\mathbb{R}^{d \times c}$.
The total scatter matrix is defined as
$$
S_{T}=\sum_{j=1}^{n}({\bf x}_{j}-{\bf \mu})({\bf x}_{j}-{\bf \mu})^{T},
$$
moreover, it is known that (Fukunaga, 1990)
$$
S_{T}=S_{W}+S_{B}.
$$
Without loss of generality, we assume that the $n$ training vectors are linear independent. Consequently, the ranks of the matrices $S_{T},~S_{B}$ and $S_{W}$ are $n-1,~c-1$ and $n-c$, respectively.

The LDA method is realized by maximizing the between-class scatter distance while minimizing the total scatter (or the within-class scatter) distance (Fukunaga 1990).
However, when the dimension of data is much larger than the number of training samples, the total scatter matrix $S_T$ (or the within scatter matrix $S_W$) will be singular, and we suffer from the small sample size problem (Fukunaga, 1990).

The null linear discriminant analysis (null LDA) method (Chen et al., 2000) is a competitive approach to overcome this difficulty. It first computes the null space of the within-class scatter matrix $S_{W}$, and then computes the principal components of the between-class scatter matrix $S_{B}$ within the null space of $S_{W}$.
In essence, the null LDA method is to find the orientation (or the transformation) matrix $W=[{\bf w}_{1}, {\bf w}_{2},\ldots,{\bf w}_{h}]\in \mathbb{R}^{d \times h}$ (of rank $h$ with $1\leq h\leq c-1$) that satisfies the following two conditions (Sharma et al., 2012)
\begin{equation}\label{eqn1.1}
S_{W}W=0, 
\end{equation}
and
\begin{equation}\label{eqn1.2}
S_{B}W\neq 0. 
\end{equation}

When $S_T$ is singular, the null LDA method solves
\begin{equation}\label{eqn1.3}
W=S_{T}^{\dag}S_{B}W, 
\end{equation}
for the orientation matrix $W$, where $S_{T}^{\dag}$ stands for the pseudo inverse (or the Moore-Penrose inverse) of $S_{T}$.
In (Sharma et al., 2012), it was shown that the equation (\ref{eqn1.3}) is a sufficient condition for the null LDA method. However, the null LDA method requires eigenvalue decomposition of $S_{T}^{\dag}S_{B}$, and the computational cost will be prohibitive when $d$ is large.
In order to release the overhead, Sharma and Paliwal (Sharma et al., 2012) propose to replace $W$ on the right-hand side of (\ref{eqn1.3}) by {\it any} random matrix $Y \in \mathbb{R}^{d \times (c-1)}$ of rank $c-1$, and make use of
\begin{equation}\label{eqn1.4}
W=S_{T}^{\dag}S_{B}Y
\end{equation}
as the orientation matrix, moreover, they present a fast implementation of null LDA method in (Sharma et al., 2012). In recent years, this method has gained wide attentions in the area of dimensionality reduction and data mining (Alvarez-Ginarte et al., 2013; Lu et al., 2013; Lyons et al., 2014; Sharma et al., 2014).

The following theorem is the main theorem of (Sharma et al., 2012). It shows that (\ref{eqn1.4}) is a sufficient condition for null LDA. Meanwhile, it is also the basis of the fast implementation of null LDA method (Sharma et al., 2012); for more details, we refer to (Sharma et al., 2012).

\begin{theorem}{\bf [Theorem 3 of (Sharma et al., 2012)]}
 If the orientation matrix $W \in \mathbb{R}^{d \times (c-1)}$ is obtained by using the relation $W=S_{T}^{\dag}S_{B}Y$ {\rm(}where $Y \in \mathbb{R}^{d \times (c-1)}$ is any random matrix of rank $c-1${\rm)}, then it satisfies the two criteria on null LDA method {\rm\big(}Eqs. {\rm(\ref{eqn1.1})} and {\rm(\ref{eqn1.2}){\rm\big)}}.
\end{theorem}
\begin{rem}
However, we find that this theorem is incomplete.
For example, let $X=[{\bf x}_1,{\bf x}_2;{\bf x}_3,{\bf x}_4]$, where $\{{\bf x}_1,{\bf x}_2\}\in\omega_1$ and $\{{\bf x}_3,{\bf x}_4\}\in\omega_2$.
Suppose that
${\bf \mu}_1=\frac{{\bf x}_1+{\bf x}_2}{2}=\hat{\bf e}~~{\rm and}~~{\bf \mu}_2=\frac{{\bf x}_3+{\bf x}_4}{2}=2\hat{\bf e}$,
where
$$
\hat{\bf e}=[1,0,1,1,\ldots,1]^T\in\mathbb{R}^d,
$$
with $d\gg n=4$. Therefore, ${\bf \mu}=({\bf x}_1+{\bf x}_2+{\bf x}_3+{\bf x}_4)/4=\frac{3}{2}\hat{\bf e}$,~
$B=[\sqrt{2}(\mu_1-\mu),\sqrt{2}(\mu_2-\mu)]=[-\frac{\sqrt{2}}{2}\hat{\bf e},\frac{\sqrt{2}}{2}\hat{\bf e}]$,
and
$$
S_B=BB^T=\hat{\bf e}\hat{\bf e}^T.
$$
Note that $rank(S_{B})=c-1=1$. In terms of Theorem 1, as $Y$ can be chosen as any random vector, we pick
$$
Y=[0,\alpha,0,\ldots,0]^T\in\mathbb{R}^d,
$$
where $\alpha$ is any positive number that satisfies $0<\alpha<1$. Then, $S_{B}Y=0$, $W=S_{T}^{\dag}S_{B}Y=0$, and
$$
S_{B}W=0,
$$
which does not satisfy the criterion {\rm (\ref{eqn1.2})}.
\end{rem}

As a result, if the random matrix is chosen arbitrarily, the orientation matrix may be rank deficient, and some discriminant information is lost.
In this paper, we revisit the fast implementation of null linear discriminant analysis method and consider how to choose the random matrix properly, such that the two criteria (\ref{eqn1.1}) and (\ref{eqn1.2}) of the null LDA method are satisfied theoretically.
We give a necessary and sufficient condition to guarantee that the orientation matrix $W$ from (\ref{eqn1.4}) is of full column rank. Moreover, the
geometric characterization of this condition is also investigated.

\section{The main result}

Since the $n$ training vectors $\{{\bf x}\}_{i=1}^n$ are linear independent, and the orientation matrix $W$ is required to be of full column rank in the null LDA method, in this paper, we focus on how to choose $Y \in \mathbb{R}^{d \times (c-1)}$ (of rank $c-1$) in (\ref{eqn1.4}), such that rank$(W)=c-1$. We follow the notations used in (Sharma et al., 2012).

Let
$$
S_{T}=U\Sigma^{2}U^{T}=
[U_{1},~U_{2}]
{\left[\begin{array}{cc}
\Sigma_{1}^{2}&0\\
0&0
\end{array}\right]}
{\left[\begin{array}{c}
U_{1}^{T}\\
U_{2}^{T}
\end{array}\right]}
$$
be the eigenvalue decomposition of $S_{T}$, where $U_{1} \in \mathbb{R}^{d \times (n-1)}$ corresponds to the range of $S_{T}$, $U_{2} \in \mathbb{R}^{d \times (d-n+1)}$ corresponds to the null space of $S_{T}$, and $\Sigma_{1} \in \mathbb{R}^{(n-1) \times (n-1)}$ is a diagonal matrix with positive diagonal elements.
From now on, we denote $G=S_{T}^{\dag}S_{B}$ for notation simplicity. By Lemma $A3$ of (Sharma et al., 2012), we have that
$$
G=S_{T}^{\dag}S_{B}
=U{\left[\begin{array}{cc}
\Sigma_{1}^{-2}&0\\
0&0
\end{array}\right]}
U^{T}S_{B}UU^{T}
=U{\left[\begin{array}{cc}
\Sigma_{1}^{-2}U_{1}^{T}S_{B}U_{1}&0\\
0&0
\end{array}\right]}U^{T},
$$
and
$$
GU=U{\left[\begin{array}{cc}
\Sigma_{1}^{-2}U_{1}^{T}S_{B}U_{1}&0\\
0&0
\end{array}\right]}.
$$
Recall that $S_{B}=BB^{T}$, thus
$$
GU=U{\left[\begin{array}{cc}
\Sigma_{1}^{-1}\Sigma_{1}^{-1}U_{1}^{T}BB^{T}U_{1}\Sigma_{1}^{-1}\Sigma_{1}&0\\
0&0
\end{array}\right]}.
$$

Let $Q=\Sigma_{1}^{-1}U_{1}^{T}B$, and $QQ^{T}=R \Lambda R^{T}$ be the eigenvalue decomposition, where $R \in \mathbb{R}^{(n-1) \times (n-1)}$ is orthonormal and $\Lambda \in \mathbb{R}^{(n-1) \times (n-1)}$ is diagonal. So we arrive at
\begin{equation}\label{eqn1.5}
G [U_{1},~U_{2}]=[U_{1},~U_{2}]
{\left[\begin{array}{cc}
\Sigma_{1}^{-1}R \Lambda R^{T}\Sigma_{1}&0\\
0&0
\end{array}\right]}.
\end{equation}
That is,
\begin{equation}\label{eqn1.6}
GU_{1}=U_{1}\Sigma_{1}^{-1}R \Lambda R^{T}\Sigma_{1},
\end{equation}
and
\begin{equation}\label{eqn1.7}
GU_{2}=0. 
\end{equation}
Moreover, it was proven in Lemma $A2$ of (Sharma et al., 2012) that $\Lambda={\left[\begin{array}{cc}
I_{c-1}&0\\
0&0
\end{array}\right]}$, where $I_{c-1}$ is the $(c-1)\times(c-1)$ identity matrix. It follows from (\ref{eqn1.6}) that
$$
GU_{1}\Sigma_{1}^{-1}R=U_{1}\Sigma_{1}^{-1}R\Lambda=U_{1}\Sigma_{1}^{-1}R
{\left[\begin{array}{cc}
I_{c-1}&0\\
0&0
\end{array}\right]}.
$$
Notice that ${\rm span}\{U_{1}\Sigma_{1}^{-1}R\}={\rm span}\{U_{1}\}$. Decompose $U_{1}\Sigma_{1}^{-1}R=[\hat{U}_{1},~\hat{U}_{2}]$, where $\hat{U}_{1}\in \mathbb{R}^{d \times (c-1)}$ is the matrix composed of the first $c-1$ columns of $U_{1}\Sigma_{1}^{-1}R$, and $\hat{U}_{2}\in \mathbb{R}^{d \times (n-c)}$, then
$$
G[\hat{U}_{1},~\hat{U}_{2}]=[\hat{U}_{1},~\hat{U}_{2}]
{\left[\begin{array}{cc}
I_{c-1}&0\\
0&0
\end{array}\right]}
=[\hat{U}_{1},~0],
$$
i.e.,
\begin{equation}\label{eqn1.8}
G\hat{U}_{1}=\hat{U}_{1}~~{\rm and}~~G\hat{U}_{2}=0.
\end{equation}
\begin{rem}
Denote $\mathcal{U}=[U_{1}\Sigma_{1}^{-1}R,~U_{2}]=[\hat{U}_{1},\hat{U}_{2},U_{2}]\in \mathbb{R}^{d \times d}$, it is seen that the columns of $\mathcal{U}$ construct a basis in $\mathbb{R}^d$, moreover, we have that
$$
{\rm rank}\big([\hat{U}_{2},U_{2}]\big)=d-c+1\gg 1.
$$
Therefore, if the $d\times (c-1)$ matrix $Y \in \emph{span}\{\hat{U}_{2},U_{2}\}$, then it follows from {\rm(\ref{eqn1.7})} and {\rm(\ref{eqn1.8})} that $W=GY=S_{T}^{\dag}S_{B}Y=0$, $S_BW=0$, and Theorem 1 fails to hold.
\end{rem}

Next, we aim to give a necessary and sufficient condition for rank$(W)=c-1$.
As the columns of $\mathcal{U}=[\hat{U}_{1},\hat{U}_{2},U_{2}]$ construct a basis of $\mathbb{R}^{d}$, for any matrix $Y  \in \mathbb{R}^{d \times (c-1)}$, there exists a matrix $[\hat{Z}_{1}^{T},\hat{Z}_{2}^{T},Z_{2}^{T}]^{T}\in \mathbb{R}^{d \times (c-1)}$, such that
\begin{equation}\label{eqn1.9}
Y=[\hat{U}_{1},\hat{U}_{2},U_{2}]
{\left[\begin{array}{c}
\hat{Z}_{1}\\
\hat{Z}_{2}\\
Z_{2}
\end{array}\right]}
=\hat{U}_{1}\hat{Z}_{1}+\hat{U}_{2}\hat{Z}_{2}+U_{2}Z_{2}. 
\end{equation}
Thus,
$$
[\hat{Z}_{1}^{T},\hat{Z}_{2}^{T},Z_{2}^{T}]^{T}=\mathcal{U}^{-1}Y,
$$
and $\hat{Z}_{1}=(\mathcal{U}^{-1}Y)(1:c-1,:)\in \mathbb{R}^{(c-1) \times (c-1)}$ is the first $c-1$ rows of $\mathcal{U}^{-1}Y$. Here $(\mathcal{U}^{-1}Y)(1:c-1,:)$ stands for the first $c-1$ rows of the matrix $\mathcal{U}^{-1}Y$.

From (\ref{eqn1.4}), (\ref{eqn1.7}), (\ref{eqn1.8}) and (\ref{eqn1.9}), we obtain
\begin{eqnarray}\label{eqn10}
W&=&S_{T}^{\dagger}S_{B}Y=S_{T}^{\dagger}S_{B}[\hat{U}_{1},\hat{U}_{2},U_{2}]
{\left[\begin{array}{c}
\hat{Z}_{1}\\
\hat{Z}_{2}\\
Z_{2}
\end{array}\right]}
=S_{T}^{\dagger}S_{B}\hat{U}_{1}\hat{Z}_{1}\nonumber\\
&=&G\hat{U}_{1}\hat{Z}_{1}=\hat{U}_{1}\hat{Z}_{1}. 
\end{eqnarray}
Since $\hat{U}_{1}$ is of full column rank, we have from (\ref{eqn10}) that $rank(W)=c-1$ if and only if $rank(\hat{Z}_{1})=c-1$, i.e., $\hat{Z}_{1}$ is nonsingular.

We are in a position to consider how to evaluate $\hat{Z}_1$ in practice. Recall that $\mathcal{U}=[\hat{U}_{1},\hat{U}_{2},U_{2}]=[U_{1}(\Sigma_{1}^{-1}R),~ U_{2}]$. Let $\Sigma_{1}^{-1}R=\hat{Q}\hat{R}$ be the QR decomposition, where $\hat{Q}\in \mathbb{R}^{(n-1) \times (n-1)}$ is an orthogonal matrix and $\hat{R} \in \mathbb{R}^{(n-1) \times (n-1)}$ is an upper triangular matrix, then
$$
\mathcal{U}=[U_1\Sigma_1^{-1}R,~U_2]=[U_1\hat{Q}\hat{R},~U_2]=[U_1\hat{Q},~U_2]{\left[\begin{array}{cc}
\hat{R}&0\\
0&I_{d-n+1}
\end{array}\right]},
$$
is the QR decomposition of $\mathcal{U}$, where $[U_1\hat{Q},~U_2]$ is orthonormal and $I_{n-d+1}$ is the $(n-d+1)\times (n-d+1)$ identity matrix. Thus,
\begin{eqnarray}\label{eqn11}
\mathcal{U}^{-1}Y&=&
{\left[\begin{array}{cc}
\hat{R}^{-1}&0\\
0&I_{n-d+1}
\end{array}\right]}
{\left[\begin{array}{c}
\hat{Q}^{T}U_{1}^{T}\\
U_{2}^{T}
\end{array}\right]}Y\nonumber\\
&=&{\left[\begin{array}{cc}
\hat{R}^{-1}&0\\
0&I_{n-d+1}
\end{array}\right]}
{\left[\begin{array}{c}
\hat{Q}^{T}U_{1}^{T}Y\\
U_{2}^{T}Y
\end{array}\right]}\nonumber\\
&=&{\left[\begin{array}{c}
\hat{R}^{-1}\hat{Q}^{T}U_{1}^{T}Y\\
U_{2}^{T}Y
\end{array}\right]}.
\end{eqnarray}
Let $\hat{R}^{-1}=
{\left[\begin{array}{c}
\hat{R}_{1}^T\\
\hat{R}_{2}^T
\end{array}\right]}$,
where $\hat{R}_{1}^T\in \mathbb{R}^{(c-1) \times (n-1)}$ is composed of the first $c-1$ rows of $\hat{R}^{-1}$, and $\hat{R}_{2}^T\in \mathbb{R}^{(n-c) \times (n-1)}$ is composed of the last $n-c$ rows of $\hat{R}^{-1}$. So we obtain from (\ref{eqn11}) that
\begin{eqnarray}
\hat{Z}_{1}&=&(\mathcal{U}^{-1}Y)(1:c-1,:)=(\hat{R}^{-1}\hat{Q}^{T}U_{1}^{T}Y)(1:c-1,:)\nonumber\\
&=&\big(U_{1}\hat{Q}\hat{R}_{1}\big)^TY. 
\end{eqnarray}

Furthermore, if $rank(\hat{Z}_{1})=c-1$, then we have from (\ref{eqn10}) that $W=S_{T}^{\dagger}S_{B}Y$ is of rank $c-1$. According to Lemma $A3$ of (Sharma et al., 2012), we have
$$
S_{T}^{\dagger}S_{B}W=(S_{T}^{\dagger}S_{B})(S_{T}^{\dagger}S_{B})Y=S_{T}^{\dagger}S_{B}Y=W,
$$
and it follows from Theorem 1 and Theorem 2 of (Sharma et al., 2012) that $W$ satisfies the null LDA criteria (\ref{eqn1.1}) and (\ref{eqn1.2}).

In summary, we have the main theorem that is a modification to Theorem 1 [Theorem 3 in (Sharma et al., 2012)].
\begin{theorem}
Let $Y \in \mathbb{R}^{d \times (c-1)}$ be a random matrix of rank $c-1$, and let
$$
\hat{Z}_{1}=\big(U_{1}\hat{Q}\hat{R}_{1}\big)^T Y \eqno(13)
$$
be the $(c-1)\times(c-1)$ matrix composed of the first $c-1$ rows of $\mathcal{U}^{-1}Y$. Then $W=S_{T}^{\dagger}S_{B}Y$ is of rank $c-1$ if and only if $\hat{Z}_{1}$ is nonsingular. Moreover, if $\hat{Z}_{1}$ is nonsigular, then $W=S_{T}^{\dagger}S_{B}Y$ satisfies the criteria of the null LDA method {\rm\big(}Eqs. {\rm (\ref{eqn1.1})} and {\rm (\ref{eqn1.2})}{\rm\big)}.
\end{theorem}

Notice that $U_{1}\hat{Q}\hat{R}_{1}$ is of full rank. Given a random matrix $Y \in \mathbb{R}^{d \times (c-1)}$, the following theorem describes the
geometric characterization of the condition for $\hat{Z}_{1}$ being nonsingular.

\begin{theorem}
Suppose that $Y \in \mathbb{R}^{d \times (c-1)}$ is of full column rank, and denote by ${\rm span}\{Y\}$ the subspace spanned by the columns of $Y$. Let $\mathcal{K}={\rm span}\{Y\}$ and $\mathcal{L}={\rm span}\{U_{1}\hat{Q}\hat{R}_{1}\}$,
then $\hat{Z}_{1}$ is nonsingular
if and only if
any nonzero vector ${\bf x}\in \mathcal{K}$ {\rm (}or ${\bf y}\in
\mathcal{L}${\rm)}, it is not orthogonal to $\mathcal{L}$ {\rm(}or
$\mathcal{K}${\rm)}.
\end{theorem}
{\bf Proof}.~~The proof is by contradiction.
On one hand, suppose that there is
a nonzero vector ${\bf x}\in \mathcal{K}$ and ${\bf x}\perp \mathcal{L}$. Then there exists a nonzero vector ${\bf z}\in
\mathbb{R}^{c-1}$, such that ${\bf x}=Y{\bf z}$. Since ${\bf x}\perp \mathcal{L}$,
we obtain
$$
0=(U_{1}\hat{Q}\hat{R}_{1})^{T}{\bf x}=\big(U_{1}\hat{Q}\hat{R}_{1}\big)^T Y{\bf z}=\hat{Z}_{1}{\bf z},
$$
and $\hat{Z}_{1}$ is singular. This
shows that, if $\hat{Z}_{1}$ is nonsingular, then
for any nonzero vector ${\bf x}\in \mathcal{K}$, it is not orthogonal to $\mathcal{L}$.
On the other hand, we assume that $\hat{Z}_{1}$ is singular. Then there is a nonzero vector ${\bf z}\in
 \mathbb{R}^{c-1}$, such that $\hat{Z}_{1}{\bf z}=\big(U_{1}\hat{Q}\hat{R}_{1}\big)^T Y{\bf z}=0$. Let ${\bf x}\equiv Y{\bf z} \in \mathcal{K}$, then ${\bf x}\neq 0$, and
it is orthogonal to $\mathcal{L}$. This implies that,
if for any nonzero vector ${\bf x}\in \mathcal{K}$, it is not orthogonal to $\mathcal{L}$, then $\hat{Z}_{1}$ is nonsingular.\quad $\Box$

\begin{rem}
Given a random matrix $Y$, Theorem 2 can be utilized to check whether $W$ is of full rank a prior in the fast implementation of the null LDA method (Sharma et al., 2012).
Indeed, it indicates that
$W$ is of rank $c-1$ if and only if $\hat{Z}_{1}=\big(U_{1}\hat{Q}\hat{R}_{1}\big)^TY$ is nonsingular.
Equivalently, Theorem 3 shows that this only happens if and only if for any nonzero vector $\bf x$ in $\mathcal{K}={\rm span}\{Y\}$, it is not orthogonal to $\mathcal{L}={\rm span}\{U_{1}\hat{Q}\hat{R}_{1}\}$, moreover, for any nonzero vector $\bf y$ in $\mathcal{L}={\rm span}\{U_{1}\hat{Q}\hat{R}_{1}\}$, it is not orthogonal to $\mathcal{K}={\rm span}\{Y\}$.

In practice, however, the case of a ``near singular" {\rm(}i.e., the smallest eigenvalue is not zero but is close to zero{\rm)} $\hat{Z}_{1}$ can occur if $Y$ is chosen arbitrarily (Golub et al., 2013). Consequently, $W$ will be near rank deficient and the two criteria of the null LDA method can not be satisfied any more.
In this situation, we suggest using another random matrix $Y$ instead.
\end{rem}

\section*{Acknowledgments}
The first author is supported by the Postgraduate Innovation Project of Jiangsu Province under grant CXLX13\_968.
The second author is
supported by the National Science Foundation of China under grant 11371176, the Natural Science Foundation of Jiangsu Province under grant BK20131126, the 333 Project of Jiangsu Province, and the Talent Introduction Program of China University of Mining and Technology.
Meanwhile, the authors would like to thank Ting-ting Xu for
helpful discussions.

\section*{References}
\bibliographystyle{model2-names}
\bibliography{refs}
\small{

 Alvarez-Ginarte Y., Montero-Cabrera L. et al., 2013. {\it  Integration of ligand and structure-based virtual screening for identification of leading anabolic steroids}, The Journal of steroid biochemistry and molecular biology, 138: 348--358.\\
 Chen L., Liao H., Ko M., Lin J. and Yu G, 2000. {\it A new LDA-based face recognition system which can solve the small sample size problem}, Pattern Recognition, 33: 1713--1726.\\
Fukunaga K., 1990. {\it Introduction to Statistical Pattern Recognition}, Academic Press Inc., Hartcourt Brace Jovanovich Publishers, SanDiego, CA 92101--4495, USA.\\
Golub G.H., Van Loan C.F., 2013. {\it Matrix Computations},
4th edition, John Hopkins University Press, Baltimore, MD.\\
Lu G., Zheng W., 2013. {\it Complexity-reduced implementations of complete and null-space-based linear discriminant analysis}, Neural Networks, 46: 165--171.\\
Lyons J., Biswas N., Sharma A. et al., 2014. {\it Protein fold recognition by alignment of amino acid residues using kernelized dynamic time warping}, Journal of theoretical biology, 354: 137--145.\\
Sharma A., Paliwal K., 2012.  {\it A new perspective to null linear discriminant analysis method and its fast implementation using random matrix multiplication with scatter matrices}, Pattern Recognition, 45: 2205--2213.\\
Sharma A., Paliwal K., Imoto S. et al., 2014. {\it  A feature selection method using improved regularized linear discriminant analysis}, Machine vision and applications, 25: 775--786.
}

\end{document}